\documentclass[12pt,leqno]{article}
\usepackage{amsmath}
\usepackage{amssymb}
\usepackage{amscd}
\def\overset#1#2{{\mathrel{\mathop {{#2}_{}}\limits^{#1}}}}
\def\underset#1#2{{\mathrel{\mathop {{}_{} {#2}}\limits_{{#1}_{}}}}}
\def\upplim_#1{\underset{#1}{\overline\lim}\;}
\def\lowlim_#1{\underset{#1}{\underline\lim}\;}
\newtheorem{cor}[equation]{Corollary}

\newtheorem{lem}[equation]{Lemma}
\newtheorem{prop}[equation]{Proposition}
\newtheorem{rmk}[equation]{\indent \rm {\it Remark}}

\newtheorem{conjecture}{\indent \rm {\it Conjecture}}
\newtheorem{thm}[equation]{Theorem}

\newcommand{\C}{{\mathbf{C}}}
\newcommand{\codim}{{\mathrm{codim}}}

\newcommand{\del}{{\partial}}

\newcommand{\kod}{{\bar\kappa}}
\newcommand{\Lie}{{\mathop{\mathrm{Lie}}}}
\newcommand{\mult}{\mathrm{mult}}
\newcommand{\N}{\mathbf{N}}

\renewcommand{\O}{{\mathcal{O}}}

\renewcommand{\P}{{\mathbf{P}}}
\newcommand{\pnc}{{\mathbf{P}^n(\mathbf{C})}}

\newcommand{\ptwoc}{{\mathbf{P}^2(\mathbf{C})}}

\newcommand{\Sing}{\mathrm{Sing}}

\newcommand{\std}{\mathrm{St}(D)}
\newcommand{\supp}{\mathrm{Supp}\,}
\newcommand{\tensor}{\otimes}
\newcommand{\Z}{\mathbf{Z}}
\newenvironment{proof}{\par{\it Proof. }}{\hfill {\it Q.E.D.}\par\vskip3pt}
\numberwithin{equation}{section}

\title{
Degeneracy of Holomorphic Curves into Algebraic Varieties
\thanks{  Research supported in part by Grant-in-Aid
   for Scientific Research (A)(1) 13304009
and (S) 17104001.}}
\author{
Junjiro Noguchi, J\"org Winkelmann and Katsutoshi Yamanoi}
\begin{document}
\setlength{\baselineskip}{18pt}
\maketitle
\begin{abstract}
Applying the Second Main Theorem of \cite{nwy04},
we deal with the algebraic degeneracy of
entire holomorphic curves $f:\C \to X$
from the complex plane
$\C$ into a complex algebraic normal variety $X$ of
positive log Kodaira dimension that
admits a finite proper morphism to a semi-abelian
variety.
We will also discuss applications to the Kobayashi hyperbolicity
problem.
\end{abstract}

\section{Introduction and main result}

Let $f: \C \to X$ be an entire holomorphic curve into a complex
algebraic variety $X$ which may be non-compact in general.
We say that $f$ is {\it algebraically degenerate}
(resp.\ {\it algebraically nondegenerate})
if there is a (resp. no) proper algebraic subset of $X$
containing the image $f(\C)$.
If $X$ is compact, we denote by $\kappa(X)$ (resp.\ $q(X)$)
the Kodaira dimension (resp.\ irregularity) of $X$;
in general, we write $\kod (X)$ (resp.\ $\bar q(X)$)
for the log Kodaira dimension
(resp.\ log irregularity) of $X$.

If $\bar q(X)> \dim X$, then every holomorphic curve $f:\C \to X$
is algebraically degenerate by log Bloch-Ochiai's Theorem
(cf.\ Theorem \ref{logbo}).
Here we would like to discuss the following problem:

{\it Problem.}  If $\bar q (X)=\dim X$ and $\kod(X)>0$,
is every holomorphic curve $f: \C \to X$ algebraically degenerate?
\smallskip

With $\bar q (X)=\dim X$ the condition of $\kod(X)> 0$
prohibits $X$ to be a semi-abelian variety, for which the Problem
clearly fails to hold.
By log Bloch-Ochiai's Theorem it is readily reduced to the case
where the quasi-Albanese map $\pi: X \to A$ is dominant.
If $\pi$ is not proper, there is a counter-example to the Problem
due to \ Dethloff-Lu \cite{dl04} Proposition~7~a) such that
$\dim X=\bar q(X)=2$ and $\kod(X)=1$.
Thus it is necessary to assume that $\pi$ is proper.
In fact it is also sufficient:

\noindent
{\bf Main Theorem }
{\it
Let $X$ be a complex algebraic variety and
let $\pi:X\to A$ be a finite morphism onto a semi-abelian variety $A$.
Let $f:\C\to X$ be an arbitrary entire holomorphic curve.
If $\kod(X)>0$, then $f$ is algebraically degenerate.

Moreover, the normalization of the Zariski closure of $f(\C)$
is a semi-abelian variety which is a finite \'etale cover
of a translate of a proper semi-abelian subvariety of $A$.
}

\begin{cor}
Let $X$ be a complex algebraic variety whose quasi-Albanese map is 
a proper map.
Assume that $\kod(X)>0$ and  $\bar q(X)\ge\dim X$.
Then every entire holomorphic curve   $f:\C\to X$ is algebraically
degenerate.
\end{cor}

We will discuss in \S4
applications for the Kobayashi hyperbolicity problem
and the complements of divisors on $\pnc$
in relation with those results which were obtained by
C. Grant \cite{g86},
Grauert \cite{g89},
Dethloff-Schumacher-Wong \cite{dswdj}, \cite{dswamj}, Dethloff-Lu
\cite{dl04} \S7,
and \cite{nwy04}.

In the surface case we provide in \S6 a result towards the strong
Green-Griffiths conjecture.
\smallskip

{\it Acknowledgement.}  We are very grateful to Professors
A. Fujiki and M. Tomari for the very helpful
suggestions on the toroidal compactifications
and Professor S.-S. Roan for encouraging discussions
on some examples of singularities.

\section{Notation and preparation}

In the present paper algebraic varieties and morphisms between them
are defined over $\C$.

For a moment let $X$ be a compact complex reduced space.
We recall some notation in the value distribution theory of
holomorphic curves (cf.\ \cite{no84}, \cite{nwy02}, \cite{nwy04}).
Let $f:\C \to X$ be an entire holomorphic curve into $X$.
Fixing a hermitian metric form $\omega$ on $X$,
we define the order function $T_f(r)$ by
$$
T_f(r)=T_f(r; \omega)=
\int_1^r \frac{dt}{t}\int_{|z|<t} f^*\omega
\quad (1\leq r<\infty).
$$

Let $V$ be a complex subspace of $X$ which may be
reducible or non-reduced.
As in \cite{nwy04} \S2 we have
\begin{enumerate}
\item
the proximity function $m_f(r; V)$ via a Weil function for $V$,
\item
the curvature form $\omega_{V,f}$ along $f:\C \to X$,
and the order function $T(r; \omega_{V,f})$ with respect to it,
\item
the counting function
$N(r;f^* V)$ and the truncated counting function $N_l(r; f^* V)$
to level $l$.
\end{enumerate}

The following F.M.T. (First Main Theorem) holds
(cf.\ \cite{nwy04} Theorem 2.9):
\begin{equation}
\label{fmt}
T(r; \omega_{V, f})=N(r; f^* V)+m_f(r; V) +O(1),
\end{equation}
provided that $f(\C) \not\subset \supp V$.

When $V$ is a Cartier divisor $D$ on $X$,
following the notation of \cite{no84},
we have the order function
$T_f(r; L(D))$ with respect to the line bundle $L(D)$ determined
by $D$ and the proximity function $m_f(r, D)$ defined
by a hermitian metric in $L(D)$.
Then we have (cf.\ \cite{nwy04} \S2)
\begin{align*}
T(r; \omega_{D,f}) &=T_f(r; L(D))+O(1), \\
m_f(r; D) &= m_f(r, D) +O(1).
\end{align*}
If $D$ is a big divisor on $X$ and $f:\C \to X$ is algebraically
nondegenerate, there is a constant $C>0$ such that
\begin{equation}
\label{big}
C^{-1}T_f(r) \leq T_f(r; L(D)) \leq C T_f(r).
\end{equation}
In this case we write
$$
T_f(r) \sim T_f(r; L(D)).
$$
Let $\lambda: X \to Y$ be a dominant rational mapping between
compact algebraic varieties $X$ and $Y$.
Let $f: \C \to X$ be algebraically nondegenerate.
Then there is a constant $C'>0$ such that
\begin{equation}
\label{order-decr}
T_{\lambda \circ f}(r) \leq C' T_{f}(r).
\end{equation}
Furthermore, if $\dim X =\dim Y$, then (cf. \cite{no84}  Lemma (6.1.5))
\begin{equation}
\label{order}
T_{\lambda \circ f}(r) \sim T_{f}(r)
\end{equation}

Suppose that $\lambda: X \to Y$ is a morphism into
a compact complex space $Y$, $W$ is a complex subspace of
$Y$ and $V \subseteq \lambda^* W$ scheme-theoretically
(i.e. the ideal sheaf defining $\lambda^*W$ should be a
subsheaf of the ideal sheaf defining $V$).
Then by \cite{nwy04} Theorem 2.9
\begin{equation}
\label{m-decr}
m_f(r; V) \leq m_{\lambda \circ f}(r; W)
\end{equation}
for a holomorphic curve $f: \C \to X$.

We recall the Main Theorem of \cite{nwy04}
in a form that we will use.

\begin{thm}
\label{smt}
{\rm (\cite{nwy04})}
Let $A$ be a semi-abelian variety and let $D$ be a reduced effective
divisor on $A$.
Then there exists an equivariant smooth compactification
$\bar A$ of $A$ such that the following inequality holds
for every algebraically nondegenerate holomorphic curve
$f:\C \to A$
\begin{align}
\label{prox}
m_{f}(r; \bar D) &= S_f(r; L(\bar D)),\\
\label{smt1}
T_f(r; L(\bar D)) &\leq N_{1}(r; f^*D)+\epsilon T_f(r; L(\bar D))
||_\epsilon, \quad \forall \epsilon>0.
\end{align}
\end{thm}

Here as usual in the Nevanlinna theory, we use the notation:
\begin{enumerate}
\item
$\log^+t=\log \max\{t, 1\}$ ($t \geq 0$),
\item
$S_f(r; L(\bar D))=O(\log^+ r)+O(\log^+ T_f(r; L(\bar D)))||$,
\item
``$||$'' (resp.\ ``$||_\epsilon$'') stands for the
inequality holds for $r>1$ outside a Borel subset of finite
Lebesgue measure (resp., where the Borel subset
depends on $\epsilon>0$).
\end{enumerate}

\begin{thm}
\label{codim2}
{\rm (\cite{nwy04})}
Let $A$ be a semi-abelian variety and let $V$ be an
algebraic subvariety of $A$. Assume that
$\codim_A V \geq 2$.
Then for every algebraically nondegenerate holomorphic curve
$f:\C \to A$
\begin{equation}
\label{codim2eq}
N(r; f^* V)  \leq \epsilon T_f(r)
||_\epsilon, \quad \forall \epsilon>0.
\end{equation}
\end{thm}

The next is called log Bloch-Ochiai's Theorem,
which we will use in a reduction.

\begin{thm}
\label{logbo}
{\rm (\cite{n77}, \cite{n81}, \cite{nw02})}
{\rm (i)}
Let $X$ be a connected compact K\"ahler manifold of dimension $n$ and
let $U$ be a complement of a proper analytic subset of $X$.
Assume that the log irregularity $\bar q(U) > n$.
Then the image of an arbitrary holomorphic curve
$f: \C \to U$ is contained in a proper analytic subset of $X$.

{\rm (ii)}
In particular, if $U$ is a semi-abelian variety (resp. quasi-torus),
the Zariski closure of the image $f(\C)$ in $X$ restricted to $U$
is a translation of a semi-abelian subvariety (resp.\ quasi-subtorus)
of $U$.
\end{thm}

We need the following result due to Kawamata.

\begin{thm} {\rm (\cite{k81} Theorem 27)}
\label{kaw2}
Let $X$ be a normal algebraic variety.
Let $\pi : X \to A$ be a finite morphism.
Then $\kod (X)\geq 0$ and there are a semi-abelian subvariety $B$ of
  $A$,
finite \'etale Galois covers $\tilde X \to X$ and
$\tilde B \to B$, and
  a normal algebraic variety $Y$ such that
\begin{enumerate}
\item
there is a finite morphism from $Y$ to the quotient $A/B$,
\item
$\tilde X$ is a fiber bundle over $Y$ with fiber $\tilde B$
and with translations by $\tilde B$ as structure group,
\item
$\kod (Y)=\dim Y=\kod (X)$.
\end{enumerate}
\end{thm}

In the special case where $\kod(X)=0$ this result takes the following
form:

\begin{thm} {\rm (\cite{k81} Theorem 26)}
\label{kaw1}
Let $X$ be a normal algebraic variety, let $A$ be a semi-abelian variety
and let $\pi: X \to A$ be a surjective finite morphism.
If $\kod(X)=0$, then $X$ is a semi-abelian variety and $\pi$
is \'etale.
\end{thm}

\section{The Compactification}
Let $X$ be a normal algebraic variety, let $A$ be a semi-abelian variety
and let $\pi:X\to A$ be a finite morphism.
We need a good compactification.

Before constructing such a compactification we remark that given
a finite morphism from a {\em normal} variety $X$ to a smooth
variety $A$ we may desingularize $X$ in order to get a generically
finite morphism between smooth varieties. Conversely, given
a generically finite morphism $p:\tilde X\to A$ between smooth
varieties, the Stein factorization gives us a normal variety $X$
together with a finite morphism from $X$ to $A$ and a proper
connected morphism from $\tilde X$ to $X$ which is a desingularization
for $X$.
\subsection{Some toric geometry}
\begin{lem}\label{finitely-many-orbits}
Let $A$ be a semi-abelian variety acting on a
(possibly singular) projective
variety $X$ with an open orbit.
Then there are only finitely many $A$-orbits in $X$.
\end{lem}
\begin{proof}
Let $\tilde X\to X$ be an equivariant desingularization.
Then  the number of $A$-orbits in $\tilde X$ is finite
(\cite{nwy04}, Lemma~3.12).
It follows that the number of $A$-orbits in $X$ is finite as well.
\end{proof}

  From toric geometry we know (see \cite{s74}, \S3):
\begin{lem}\label{aff-nbghd}
Let $Y$ be a normal toric variety with algebraic torus $T$.
Then every point in $Y$ admits a Zariski open invariant affine
neighbourhood $W$.
\end{lem}

\begin{lem}\label{extend-inv}
Let $W$ be a Stein complex space on which a reductive complex
Lie group $H$ acts. Let $V$ be a complex vector space,
let $Z\subset W$ be an invariant closed analytic subset and
let $f_0:Z\to V$ be an $H$-invariant holomorphic map.
Then $f_0$ extends to an $H$-invariant holomorphic map $f:W\to V$.
\end{lem}
\begin{proof}
Let $f_1:W\to V$ be an arbitrary (i.e. not necessarily
invariant)  holomorphic extension (which exists, because $W$
is Stein). Let $K$ be a maximal compact subgroup of $T$ with
normalized Haar measure $d\mu$.
Then we can define $f$ by
\[
f(x)=\int_K f_1(k\cdot x)d\mu(k)
\]
\end{proof}

\begin{lem}\label{lem-aff}
Let $W$ be an irreducible
affine variety on which a reductive commutative
complex Lie group
$T=(\C^*)^g$ acts with finitely many orbits.
Let
$S_i$, $i=1,2$ be two $T$-orbits.
Then either $S_1=S_2$ or the isotropy Lie algebras differ.
\end{lem}
(Note that two points $p$, $q$ in the same $T$-orbit $S$ have
the same isotropy group, because $T$ is commutative.
Therefore it makes sense to talk about the ``isotropy Lie algebra
of an orbit''.)
\begin{proof}
Let $S_1,S_2$ be two orbits with the same isotropy Lie algebra.
Note that this implies $\dim(S_1)=\dim(S_2)$.
Let $Y$ denote the closure of $S_1\cup S_2$ in $W$.
Let $H$ denote the connected component of the isotropy group
for $S_1$. By our assumption $H$ acts trivially on $Y$.
Now we choose a holomorphic embedding
$i:Y\hookrightarrow \C^N$. Since $H$ acts trivially on $Y$,
this map $i$ is $H$-invariant. Thus we can extend it to
an $H$-invariant holomorphic map $F:W\to\C^N$
(using Lemma~\ref{extend-inv}).
We observe that
\[
\dim(F(Y))=\dim(i(S_1))=\dim T-\dim H.
\]
Since $F$ is $H$-invariant,
\[
\dim F(\Omega )\le \dim T-\dim H
\]
for every $T$-orbit $\Omega\subset W$.
Recall that we assumed that the number of $T$-orbits is finite.
Therefore the above inequality implies $\dim F(W)\le\dim T-\dim H$.
It follows that $\dim F(Y)=\dim F(W)$. Since $W$ is irreducible
and $Y\subset W$, it follows that $F(Y)=i(Y)\simeq Y$
is irreducible. Therefore $S_1=S_2$.
\end{proof}

\begin{lem}\label{iso-lin}
Let $A$ be a semi-abelian variety
acting effectively on an algebraic variety $X$.
Let $L$ be the maximal connected linear subgroup of $A$ and $x\in X$.
Then $L$ contains the connected component of the isotropy group
$A_x=\{a\in A:a\cdot x=x\}$.
\end{lem}
\begin{proof}
Let $\O_x$ be the local ring at $x$ and let $m_x$ be its maximal ideal.
Then $A_x$ acts linearly on each vector space $\O_x/m_x^k$.
Because the $A$-action is supposed to be effective, these actions
can not be all trivial.
Hence $A_x$ is linear and its connected component is contained in $L$.
\end{proof}

\begin{lem}\label{alb-sing}
Let $A$ be a semi-abelian variety
acting on a normal algebraic variety $X$ with an open orbit.
Let $\pi:\tilde X\to X$ be an equivariant desingularization
and $a:\tilde X\to T$ the Albanese map.
Then there exists an equivariant morphism $f:X\to T$
such that $a=f\circ\pi$.
\end{lem}
\begin{proof}
First we recall that there are only finitely many $A$-orbits in
$X$ and $\tilde X$ (Lemma~\ref{finitely-many-orbits}).
Let $L$ be the maximal linear connected subgroup of $A$.
Let $p\in X$, let $\Omega=A(p)$ and let
$\tilde\Omega$ be an irreducible component of $\pi^{-1}(\Omega)$.
Since $\tilde\Omega$ is invariant and irreducible and contains only
finitely many $A$-orbits, it contains a dense $A$-orbit.
Furthermore observe that the fibers of $\pi$ are irreducible,
because $X$ is normal. As a consequence we obtain that the fiber
$\pi^{-1}(p)$ equals the closure of an orbit of the connected
component of the isotropy
group $A_p^0$ acting on $\tilde\Omega$.
Since  $A_p^0\subset L$ (see Lemma~\ref{iso-lin}), we see that all fibers of $\pi$ are contained
in closures of $L$-orbits in $\tilde X$.
Being linear, $L$ acts trivially on the abelian variety $T$.
Therefore $a:\tilde X\to T$ is constant along the fibers of $\pi$.
Since $X$ is normal this implies that $a$ fibers through $\pi$.
\end{proof}

\begin{prop}\label{prop-toric}
Let $A$ be a semi-abelian variety acting on a normal
algebraic variety $X$ with an open orbit.
Then every point $x\in X$ admits a Zariski open neighbourhood
$W$ of $x$
in $X$ such that the following property holds:

``Two points $y,z\in W$ are contained in the same $A$-orbit
if and only if they have the same isotropy Lie algebra with respect
to the $A$-action.''
\end{prop}
\begin{proof}
We may assume that the $A$-action on $X$ is effective.
Let $L$ be the maximal connected linear subgroup of $A$
and $T=A/L$. Due to Lemma~\ref{alb-sing}
there is a surjective equivariant morphism
$\alpha:X\to T$ (which is the Albanese of a desingularization
of $X$).

Now $A$-invariant subsets of $X$ correspond
to $L$-invariant subsets of a fiber $F$ of this
morphism $\alpha:X\to T$.
Moreover, due to Lemma~\ref{iso-lin} every isotropy Lie algebra
for the $A$-action is contained in the Lie algebra $Lie(L)$ and
therefore coincides with the isotropy Lie algebra  for the $L$-action.
Therefore there is no loss in generality in assuming that $L=A$.
Then $X$ is a toric variety and every point $x\in X$ admits
an invariant affine neighbourhood (Lemma~\ref{aff-nbghd}).
Due to Lemma~\ref{finitely-many-orbits} there are only finitely
many orbits.
Hence the statement follows from Lemma~\ref{lem-aff}.
\end{proof}

\begin{lem}\label{comm-action}
Let $G$ be a commutative algebraic group acting with
a dense connected open orbit $\Omega$ on a variety $X$ and
acting transitively on a
variety $Y$. Let $\pi:X\to Y$ be a surjective equivariant
morphism.
Then for every $y\in Y$ and every $G$-orbit $Z\subset X$ the
intersection with the fiber $Z\cap\pi^{-1}(y)$ is connected.
\end{lem}
\begin{proof}
Fix $y\in Y$ and define $F=\pi^{-1}(y)$. Let $H$ denote the isotropy group
$H=\{g\in G:g(y)=y\}$. There is a one-to-one correspondence between
$G$-orbits in $X$ and $H$-orbits in $F$.
$W=\Omega\cap F$ is a dense open $H$-orbit in $F$.
Choose $p\in W$ and let $I=\{h\in H:h(p)=p\}$. Because $H$ is commutative,
the group $I$ acts trivially on $W$. Since $W$ is dense in $F$, the action
of $I$ on $F$ is trivial. Let $H^0$ denote the connected component of $H$.
Then $H^0\cdot I=H$, because $W\simeq H/I$ is connected.
Thus the fact that $I$ acts trivially on $F$ implies that the $H$-orbits
on $F$ coincide with the $H^0$-orbits. In particular they are
connected. Thus $Z\cap F$ which is an $H$-orbit
must be connected.
\end{proof}

\begin{lem}\label{RS}
Let $\pi:X\to Y$ be a finite morphism from a normal variety $X$
with singular locus $S$ onto a smooth variety $A$.
Let $R$ denote the ramification divisor of the restriction of $\pi$
to $X\setminus S$. Then $S$ is contained in the closure of $R$.
\end{lem}

\begin{proof}
Let $p\in X\setminus\bar R$. Then there are small
 connected open Stein neighbourhoods $V$ of $p$ in $X\setminus\bar R$
and $W$
of $\pi(p)$ in $A$ such that $\pi$ restricts to a finite morphism
from $V$ to $W$. We may assume that $W$ is simply-connected.
Because $X$ is normal, the codimension of
$\pi(S)\cap W$ is at least two. Hence $W\setminus\pi(S)$ is simply-connected.
It follows that there is a section $\sigma:W\setminus\pi(S)\to V$
which (again because of $\codim\pi(S)\ge 2$) extends to all of $W$.
This yields a biholomorphic map between $V$ and $W$. Since $W\subset A$
is smooth, we deduce $p\not\in S$.
\end{proof}

\subsection{Simple compactification}
\begin{prop}\label{simple-cpt}
Let $\pi:X\to A$ be a finite surjective morphism between normal
varieties and let $A\hookrightarrow \bar A$ be
a normal compactification of $A$.
Then there exists a unique normal compactification $X\hookrightarrow
\bar X$ such that $\bar X$ is normal and $\pi$ extends to a finite
morphism $\bar\pi:\bar X\to\bar A$.
\end{prop}
\begin{proof}
Let $\Gamma\subset X\times A$ be the graph of $\pi$.
Choose a compactification $X\hookrightarrow\tilde X$ and
let $\hat\Gamma$ be the closure of $\Gamma$ in $\bar A\times\tilde X$.
Then $\bar X$ is obtained by first normalizing $\hat\Gamma$ and then
taking the Stein factorization of the projection onto $\bar A$.
It is easy to deduce unicity from the assumption of $X$ being normal.
\end{proof}

\subsection{A better compactification}
\begin{prop}\label{good-cptf}
Let $\pi$ be a finite morphism from a normal algebraic variety $X$
to a semi-abelian variety $A$.
Let $R$ denote the set of all non-singular points $p\in X$
at which $\pi$ is ramified.
Let $A\hookrightarrow \bar A$ be a smooth equivariant compactification
and let $\omega$ be a log volume form on $\bar A$.
Then there exists a compactification $X\hookrightarrow\bar X$
and a proper morphism $\bar \pi:\bar X\to\bar A$
such that:
\begin{enumerate}
\item
$\bar X\setminus\bar R$ is smooth and
$\bar X\setminus (X\cup\bar R)$ is
an s.n.c. (=``simple normal crossing'')
divisor
(where $\bar R$ denotes the closure of $R$ in $\bar X$),
\item
$\bar\pi|_X=\pi$,
\item
Let $\tilde\omega=
\bar\tau^*\omega\in\Omega^d(\bar X;\log \del X)$.
Then $\tilde\omega$ has poles along all divisorial components of
$\bar X\setminus(X\cup\bar R)$.
\end{enumerate}
\end{prop}

\begin{proof}
Due to Lemma~\ref{RS} the singular locus $S$ of $X$ is contained
in the closure of $R$.
We let $D$ be a divisor on $A$ containing $\pi(R)$.
Then $D$ contains $\pi(S)$, too.

(0) To prove the assertion we use the following strategy:
\begin{itemize}
\item
We define a class of ``admissible compactifications'' of $X$.
\item
We show, using Lemma~\ref{simple-cpt}, that there exists an
admissible compactification.
\item
For each admissible compactification we define an
 ``indicator function''
$\zeta:\partial A\to \N$ which measure the presence of singularities
outside of $\bar R$.
\item
Using the theory of toroidal embeddings, we show that we can blow up
admissible compactification in such a way that we stay inside
the category of admissible compactification, but decrease the
indicator function.
\item
We verify that after finitely many steps the indicator function
vanishes and that then we have found a compactification as desired.
\end{itemize}

(1)
A compactification
$X\hookrightarrow X'$ is
``admissible'' if the following properties are satisfied:

\begin{enumerate}
\item
$X'$ is normal,
\item
the projection map $\pi:X\to A$ extends to a proper holomorphic map
$\pi':X'\to \bar A$ with ${(\pi')}^{-1}(A)=X$.
\item
For each point $p\in\ X'\setminus\bar R$ there is an open neighbourhood
$U$ of $\pi(p)=q$ in $\bar A$ such that the connected component
$\Omega$ of ${(\pi')}^{-1}(U)\setminus\bar R$
which contains $p$ admits a biholomorphic
map $\psi:\Omega\to W$ into an open subset $W$ of a toric variety $Z$.
\item Let $G$ denote the algebraic torus $(\C^*)^g$ acting on the
toric variety $Z$.
Then map $\pi\circ\psi^{-1}:W\to \bar A$ is (locally) equivariant
for some holomorphic Lie group homomorphism with discrete fibers
from  $G$ to $A$.
\end{enumerate}

Condition (iii) could be rephrased by saying that
$X\hookrightarrow X'$
should be locally a ``toroidal embedding'' in the sense of
\cite{kkms73} except at $\bar R$.

(2)
Using Lemma~\ref{simple-cpt},
we obtain  a normal compactification
$X\hookrightarrow X_1$ such that $\pi:X\to A$ extends to a finite
morphism $\pi_0:X_1\to\bar A$.
Let $p\in X'\setminus\bar R$, $q={(\pi')}(p)$.
Let $U$ be an open neighbourhood
of $p$ in $X'\setminus\bar R$ and define
$\Omega={(\pi')}(U)$.
Since $\partial A$ is an s.n.c. divisor,
we can shrink these open neighborhoods to obtain local coordinates
$z_i$ on $\Omega$  such that
$z_i(q)=0$ and such that $\partial A$ is the zero locus of
$z_1 \cdots z_k$ for some $k$.
By shrinking $\Omega$ we may assume that $\Omega$ is biholomorphic
to a polydisc with $\Omega\cap A\simeq(\Delta^*)^k\times\Delta^{n-k}$.
Let $i:\Delta^*\to\C^*$ be the standard injection and $j:\Delta\to\C^*$
any open embedding (e.g.~$j(z)=z+2$).
We obtain $\xi=(i^k,j^{n-k}):\Omega\cap A\hookrightarrow(\C^*)^n=G$.
Now $\pi_1(G)\simeq\Z^n$ and $\xi$ induces an embedding of
$\pi_1(\Omega\cap A)\simeq\Z^k$ into $\Z^n\simeq\pi_1(G)\simeq\Z^n$
as $\xi_*\left(\pi_1(\Omega\cap A)\right)=\Z^k\times\{0\}^{n-k}$.
Thus for each subgroup $\Gamma$
of fixed finite index $d$ of $\pi_1(\Omega\cap A)$
we can choose a corresponding subgroup of the same
finite index in $\pi_1(G)\simeq\Z^n$, namely
$\xi_*(\Gamma)\times\Z^{n-k}$.
Therefore the unramified covering $\pi^{-1}(\Omega\cap A)\to\Omega\cap A$
extends via $\xi$ to an unramified covering $G_1\to G$ of $G$.
Now $G_1$ is again an algebraic group and as algebraic group
isomorphic to $(\C^*)^n$.
Since the procedure of normal  compactification as
in Proposition~\ref{simple-cpt}
is canonical, we can embed $\pi^{-1}(\Omega)\to\Omega$
into a finite morphism $\bar G_1\to\bar G=(\P_1)^n$. Now $\bar G_1$
is a toric variety since the $G_1$ action on itself extends to the
boundary $\partial G_1\simeq\partial G$ via $G_1\to G$.

This proves that $X\hookrightarrow X_1$ is admissible.

(3) We are looking for admissible compactifications which are
smooth outside the closure of $R$.
Thus given an admissible compactification $X\hookrightarrow\bar X$
with $\pi:\bar X\to\bar A$
we define our indicator function $\zeta:\partial A\to \N$ as follows:
$\zeta(p)$ denotes the number of connected components of the fiber
$\pi^{-1}(p)$ which intersect $\Sing(\bar X)\setminus\bar R$.
Evidently this function $\zeta$ vanishes iff $\bar X\setminus\bar R$
is smooth.

(4) We suppose given an admissible compactification $\bar X$ with
indicator function $\zeta$.
Since the level sets
$\{z\in\partial A:\zeta(z)=c\}$ ($c\in\N$)
of $\zeta$ are constructible sets,
it makes sense to define $d_\zeta=\dim \{z:\zeta(z)\ne 0\}$
and to define a number $n_\zeta$ which is the generic value of $\zeta$
on $\{z:\zeta(z)\ne 0\}$.
We choose a generic point $p\in\{z\in\partial A:\zeta(z)\ne 0\}$
and a point $q\in\pi^{-1}(p)\cap\Sing(\bar X)\setminus\bar R$.
Since the compactification is admissible, there is an isomorphism
$\phi:U\to W$ where $U$ is a neighbourhood of the connected component
of $\pi^{-1}(p)$ containing $q$ and $W$ is an open neighbourhood
of $p$ in a toric variety $Z$.
By the theory of toroidal embeddings
(see e.g. \cite{kkms73}) there exists an equivariant
desingularization $\tilde Z\to Z$, and thus a desingularization of $U$.
The problem is to extend this blow-up of $U$ to a blow-up of $\bar X$.
Because the blow-up of $Z$ is equivariant and the number of $G$-orbits
in $Z$ is finite, it is given by a blow-up of invariant strata
(=closures of $G$-orbits in $Z$). Now $\bar X$ admits a natural
stratification induced by the $A$-action on $\bar A$. In fact,
there is a local $A$-action on $\bar X\setminus\bar R$ which gives
this stratification.
In order to extend the blow-up, we need to extend its center,
and this means: Given a closed invariant subvariety $Q\subset Z$
we need to show that either $M\cap U=\phi^{-1}(Q)$
 or $M\cap U=\emptyset$ for each stratum $M$ of $\bar X$.%
\footnote{In the language of \cite{kkms73}, here we prove that
the embedding of $X\setminus R$ into $\bar X\setminus\bar R$
is a toroidal embedding {\em ``without self-intersection''.}}
We prove this indirectly. So let us assume that this property fails.
Then $\phi(M\cap U)$ intersects several $G$-orbits some
of which are contained in $Q$ and some of which are not.
Since $M$ is one stratum, $\pi(M)$ is one $A$-orbit in $\bar A$
which implies that for each point in $M$ we have the same isotropy
Lie algebra.
If we have chosen $U$ sufficiently small,
we can conclude from Proposition \ref{prop-toric}, that
$\pi(M\cap U)$ is connected.
Let $N=\pi^{-1}(\pi(M))$. Using Lemma~\ref{comm-action}
we may deduce that $M\cap U$ is connected, as desired.

(5)  We recall that we introduced an indicator function $\zeta$
with associated numbers $d_\zeta$ and $n_\zeta$.
For a generic point $p\in\Sing(\bar X)\setminus\bar R$
the above considerations show that there is an appropriate blow-up
yielding an other admissible compactification
which is smooth around (the preimage of) $p$.
By the definition of $d_\zeta$ and $n_\zeta$ this means that
given an admissible compactification we can always blow-up
$\bar X$ so that
either $d_\zeta$ decreases or $n_\zeta$ decreases while
$d_\zeta$ is kept fixed.
Thus we can strictly decrease the value of $(d_\zeta,n_\zeta)\in\N^2$
where $\N^2$ is endowed with the lexicographic
order. It follows that $(d_\zeta,n_\zeta)=(0,0)$ after finitely many steps.
But $(d_\zeta,n_\zeta)=(0,0)$ implies the vanishing of $\zeta:\partial A\to\N$
which in turn implies that $\bar X$ is smooth
outside $\bar R$.

Thus we have established: There exists an admissible
compactification $X\hookrightarrow X'$ such that
$X'$ is smooth outside the closure of $R$.

(6) By the definition of admissibility $X'\setminus\bar R$
is locally a toroidal embedding.
This implies that $\partial X'\setminus\bar R$ is an s.n.c.
divisor in $X'\setminus\bar R$.
Furthermore: The logarithmic tangent bundle on a smooth
toric variety is trivial.
Hence conditions (iii) and (iv) of ``admissibility''
imply the statement about the poles of $\tilde\omega$.
\end{proof}

\subsection{The best compactification}
\begin{thm}\label{thm-cptf}
\label{des-cpt}
Let $\pi: X \to A$ be a finite morphism from  a normal variety $X$
onto a semi-abelian variety $A$.
Let $\bar A$ be a smooth equivariant compactification of $A$.
Let $D$ denote the critical locus, i.e. the closure
of the set of all $\pi(z)$
where $z$ is a smooth point and
 $d\pi: T_zX\to T_{\pi(z)}A$ fails to have full rank.
Then there exist
\begin{enumerate}
\item[{\rm (a)}]
A desingularization $\tau:\tilde X\to X$ and a smooth compactification
$j:\tilde X\hookrightarrow\hat X$  such that
the boundary divisor $\partial\hat X=\hat X\setminus j(\tilde X)$ has
only simple normal crossings,
\item[{\rm (b)}]
a proper holomorphic map $\psi:\hat X\to\bar A$
such that $\psi\circ j=\pi\circ\phi$ with $\psi^{-1}(A)=\tilde X$.
\item[{\rm (c)}]
an effective divisor $\Theta$ on $\hat X$,
\item[{\rm (d)}]
a subvariety $\hat S\subset\hat X$
\end{enumerate}
such that
\begin{enumerate}
\item
$\Theta$ is linearly equivalent to the log canonical divisor
$K_{\hat X}+  \del \hat X$ of $\tilde X $,
\item
$\psi(\supp \Theta) \subset \bar D$,
\item
the image $\psi(\hat S)$ has at least codimension two in $A$,
\item
for every holomorphic curve $f:\Delta \to A$
from a disk in $\C$ with lifting $F:\Delta \to\tilde X$
and for $z\in F^{-1}(\supp \Theta\setminus \hat S)$
we have the following inequality of multiplicities:
\end{enumerate}
\begin{equation}
\label{mult}
\mult_z F^* \Theta \le \mult_z f^* D  -1.
\end{equation}
\end{thm}

\begin{proof}
Let $\Sing(X)$ be the singular locus of $X$ and
let $R$ be the ramification divisor of $\pi$ restricted
to $X\setminus\Sing(X)$.
We apply Proposition \ref{good-cptf} and obtain
a first compactification $X'$ of $X$ and an extension
$\pi':X'\to\bar A$ of $\pi:X\to A$.
We recall that $X'\setminus\bar R$ is smooth (assertion (i) of
Proposition \ref{good-cptf}).
Thus we can use Hironaka desingularization to desingularize $X'$
without changing $X'\setminus\bar R$.
We obtain a desingularization $\tau:\hat X\to X'$ which restricts
to a desingularization $\tau_0:\tilde X\to X$.

Set
$$
\Theta=(\psi^*\omega)+\psi^{-1}\del A .
$$
Then the properties (i) and (ii) are satisfied.

Let $S_1$ denote the image of the set of singular points of $X$ by
$\pi$,
let $S_2$ be the set of singularities of $D$, and let $S_3$ be the image
of the set of singularities of $\supp \Theta|_{\tilde X}$ by $\psi$.
Then $\codim\, S_j \geq 2$ $(1 \leq j \leq 3)$.
Set $S=S_1 \cup S_2 \cup S_3$ and $\hat S=\psi^{-1}(\bar S)$.
Then $\codim\, \hat S \geq 2$ and (iii) is satisfied.

Let $p \in \supp \Theta |_{\tilde X} \setminus \hat S$.
Then $q=\pi(p) \in D \setminus S$.
Because of the construction there exist local coordinate systems
$(x_1, \ldots, x_n)$ about $p$ and $(y_1, \ldots, y_n)$ about $q$
such that locally
\begin{equation*}
\begin{array}{ll}
y_1 =x_1^k \quad (k \in \N, k \geq 2), & y_j=x_j \: (2 \leq j \leq n),\\
  \\
\Theta =(k-1)\{x_1=0\},  & D=\{y_1=0\}.
\end{array}
\end{equation*}
Set $F=(F_1, \ldots, F_n)$ and $f=(f_1, \ldots, f_n)$.
Let $\nu$ denote the multiplicity of zero of $F_1$ at $z$.
Note that $\mult_z F^*\Theta=(k-1)\nu$ and $\mult_z f^*D=k\nu$.
Thus $(k-1)\nu \leq k\nu -1$ and so \eqref{mult} holds.
\end{proof}

{\it Remark.}
Let the notation be as in Proposition \ref{des-cpt}.
Let $\std$ be the stabilizer of $D$ in $A$ defined by the identity
component of $\{a \in A; a+D=D\}$.
Assume that $X$ is of log general type.
We claim
\begin{equation}
\label{stab}
\std=\{0\}.
\end{equation}
Suppose that it is not the case.
Then $ A \setminus D$ is also $\std$-invariant.
Then, setting $U=X \setminus \pi^{-1}D$, we have a sequence of
surjective morphisms
$$
U \:\: \overset{\pi|_{U_{\empty}}}{\to}  \:\:  A\setminus D \to
 (A/\std)\setminus (D/\std)
$$
induced from the restriction of $\pi$ and the quotient map.
Let $\nu: U  \to (A/\std)\setminus (D/\std)$ be the composed morphism.
Then every connected component $Z$ of the fibers of $\nu$ is
a finite \'etale cover over $\std$.
Since $\kod(Z)=0$, $\kod(U)\leq \dim U -1$ (Iitaka \cite{iitaka} Theorem 4).
Then we have a contradiction, $\kod(X)\leq \dim X -1$.
By \eqref{stab} $\bar D$ is big on $\bar A$
(\cite{nwy04} Proposition 3.9).

\section{Proof of the Main Theorem}

The following is an essential case.
\begin{thm}
\label{gen}
Let $X$ be a normal variety of log general type and
$\pi:X\to A$ be a finite morphism onto a semi-abelian variety.
Then every holomorphic curve $g:\C \to X$ is algebraically degenerate.
\end{thm}
\begin{proof}
We use the desingularization and
compactification obtained by Theorem~\ref{thm-cptf}
and follow the notation there.

Assume that there exists an algebraically nondegenerate
holomorphic curve $g:\C\to X$.
Since $g$ is algebraically nondegenerate and $\tilde X \to X$ is birational,
we can lift $g$ to an algebraically nondegenerate entire curve
$F: \C \to \tilde X$. Set $f=\pi\circ g$.
We are going to deduce a contradictory estimate for the order
function $T_F(r)$ of $F$.

Let $\bar A$ be an equivariant compactification of $A$
such that $\bar D$ is in general position;
that is, $\bar D$ contains no $A$-orbit in $\bar A$.

Since $\bar D$ is in general position in $\bar A$,
it follows from \eqref{m-decr}, \eqref{prox} and
Theorem \ref{thm-cptf} (ii), (iii) that
\begin{equation}
\label{Theta1}
m_F(r; \Theta)\leq m_f(r; \bar D)=S_f(r).
\end{equation}
Since $\Theta$ is big, one infers from \eqref{order} that
\begin{equation}
\label{order-fF}
T_f(r) \sim T_F(r).
\end{equation}
Combining this with \eqref{Theta1}, one gets
\begin{equation}
\label{Theta2}
m_F(r; \Theta)=S_F(r).
\end{equation}

Theorem~\ref{thm-cptf} (iv)
implies
\begin{equation}
\label{N-Theta1}
N(r; F^*\Theta) \leq N(r; F^* \hat S)+N(r; f^* D)-N_1(r; f^*D).
\end{equation}
Now $\psi(\hat S)$ is of codimension at least two in $A$.
Therefore we can infer from Theorem~\ref{codim2}
that
\begin{equation}
\label{hatS}
N(r; F^* \hat S) \leq N(r; f^* (\psi_* \hat S))
 \le  \epsilon T_f(r)||_\epsilon.
\end{equation}

By virtue of \eqref{smt1} we have
\begin{align*}
N(r; f^* D)-N_1(r; f^*D) &\leq T_f(r; L(\bar D))-N_1(r; f^*D)\\
&\leq \epsilon T_f(r)||_\epsilon, \qquad \forall \epsilon>0.
\end{align*}
The combination of this with \eqref{order-fF} yields
\begin{equation}
\label{ND}
N(r; f^* D)-N_1(r; f^*D) \leq \epsilon T_F(r)||_\epsilon,
 \quad \forall \epsilon>0.
\end{equation}

Now one infers from \eqref{N-Theta1}--\eqref{ND} that
\begin{equation}
\label{N-Theta2}
N(r; F^*\Theta) \leq \epsilon T_F(r)||_\epsilon,
 \quad \forall \epsilon>0.
\end{equation}
Note that $T_F(r)\sim T_F(r; L(\Theta)) \to \infty$ as $r \to \infty$.
Henceforth, F.M.T. \eqref{fmt}, \eqref{Theta2} and \eqref{N-Theta2}
lead to a contradiction:
$$
T_F(r)\leq \epsilon T_F(r)||_\epsilon,
 \quad \forall \epsilon>0.
$$
\end{proof}

{\it Remark.}
The case where $X$ is compact
was proved by Yamanoi \cite{ya04} Corollary 3.1.14.
\medskip

{\it Proof of the Main Theorem.}  Let the notation be as in
the Main Theorem.
Assume that $\bar \kappa(X)>0$ and $f: \C \to X$ is algebraically
nondegenerate.
By lifting $f$ to the normalization of $X$ we may assume further
that $X$ is normal.

We use Kawamata's Theorem \ref{kaw2} and the notation there.
Since $\tilde X\to X$ is \'etale, we can lift $f$
to a holomorphic curve $\tilde f:\C\to\tilde X$.
By Theorem \ref{gen} the composed map of $\tilde f$ with $\tilde X \to Y$
must be algebraically degenerate.
This implies that $f$ itself is algebraically degenerate, because
$\dim Y =\kod(X)>0$ due to our assumption; this is a contradiction.

Now let $Z$ be the normalization of the Zariski closure of
$f(\C)$ in $A$.
Then $\bar \kappa (Z)=0$ by what we have just proved, and
the last statement follows from Kawamata's Theorem~\ref{kaw1}.
{\it Q.E.D.}

\section{Applications}

Here we give several applications.
As a direct consequence of the Main Theorem we have the next:

\begin{thm}\label{hyp2}
Let $A$ be a semi-abelian variety with smooth equivariant
compactification $A\hookrightarrow\bar A$ and let $\bar X$
be a projective variety with a finite morphism $\pi:\bar X\to\bar A$.
Then $X=\pi^{-1}(A)$ is Kobayashi hyperbolic and hyperbolically embedded
into $\bar X$ unless there exists a semi-abelian subvariety $B\subset A$,
a positive-dimensional orbit $B(p)\subset\bar A$, an \'etale cover
$\rho:C\to B(p)$ from a semi-abelian variety $C$
and a morphism $\tau:C\to \bar X$ such that $\rho=\pi\circ\tau$.
\end{thm}

\begin{proof}
By Brody-Green's Theorem either $X$ is hyperbolic and
 hyperbolically embedded into $\bar X$
or there is a non-constant holomorphic map from $\C$ into one
of the strata of the natural stratification on $\bar X$,
the one induced by the stratification on $\bar A$ which is given by
the $A$-orbits.
Using this, the statement follows from the Main Theorem.
\end{proof}

\begin{thm}
Let $X$ be a normal algebraic variety which admits a finite
morphism $\pi$ onto a simple abelian variety $A$.
Then either $X$ is hyperbolic or $X$ itself is an abelian variety.
\end{thm}

\begin{proof}
Since $A$ is assumed to be simple, there does not exist any
non-trivial (semi-) abelian subvariety.
Hence every holomorphic map $f:\C\to X$ must be constant unless
$\pi:X\to A$ is an \'etale cover by the Main Theorem.
If every holomorphic map $f:\C\to X$ is constant, then $X$ is hyperbolic
by Brody's theorem. If $\pi$ is \'etale, then $X$ is an abelian variety, too.
\end{proof}

{\it Remark.}
(i) The case of $\dim X=\kappa (X)=2$ was proved
by C. Grant \cite{g86}.

(ii) If $A$ is not simple, it is not sufficient to assume
the Kobayashi hyperbolicity of the ramification locus $R$
of $\pi$ to obtain the Kobayashi hyperbolicity of $X$,
as shown by the next proposition.

\begin{prop}
\label{count}
There exists an abelian surface $A$ with a smooth ample hyperbolic
curve  $D \subset A$ and a
smooth projective surface $X$ with a finite covering
$\pi:X\to A$ with ramification locus $R$
such that $R=\pi^{-1}(D)$
and $X$ is not hyperbolic.
\end{prop}

\begin{proof}
Let $E$ be an elliptic curve with line bundles $H$ and $H'$
of degree
$2$ and $3$ respectively. Then $\phi_H:E\to\P_1$ while
$\phi_{H'}$ embeds $E$ into $\P_2$ as a cubic curve $C$.
Let $X=\P_1\times C$. Furthermore let $Z$ be the union of all
$\{p\}\times C$ for points $p\in\P_1$ over which $\phi_H:E\to\P_1$
is ramified.
Choose an even number $d$.
Then using Bertini's theorem we obtain a hypersurface
$L\subset X$ such that $L$ is ample, smooth, with $L\cap Z$ being
smooth too, and such that the bidegree is
$(1,d)$.
Let $\tau:E\times E\to X=\P_1\times C$ be given as
$\tau(x_1,x_2)=(\phi_H(x_1),\phi_{H'}(x_2))$.
Then $D=\tau^*L$ is a smooth ample divisor
(smoothness of $D$ can be deduced from the conditions that both $L$ and
$L\cap Z$ are smooth)
and by construction $L(D)=p_1^*L(H)\tensor p_2^*L(dH')$
where $p_i:E\times E\to E$ are the respective projections.
Since $D$ is smooth and ample in $E\times E$, it is a curve of
genus larger than one and therefore hyperbolic.
By construction there is a divisor $D_0$ such that $2D_0$ is
linearly equivalent to $D$. Now, by the usual cyclic covering
method there is a surface $X$ with a two-to-one covering
$\pi:X\to A=E\times E$ which is precisely ramified over $D$.
More precisely, by taking squares fiber-wise there is a morphism
  from (the total space of) $L(D_0)$ to $L(D)$. We let $\sigma$
denote a section of $L(D)$ whose zero-divisor is $D$ and
define $X=\{(x,t)\in L(D_0):t^2=\sigma(x)\}$.
Now $d(t^2-s(x))=2tdt-ds$ in local coordinates; hence $X$ is
smooth if $D$ is smooth.

We claim that $X$ is not hyperbolic.
  For each $q\in E$ let $E_q=E\times\{q\}\subset A$.
Then $D\cap E_q$ is a divisor of degree $2$ which is linearly
equivalent to $H$. If in some neighborhood of $q$ there
are holomorphic functions $a,b$ with values in $E$ such that
$\{a(p),b(p)\}=D\cap E_p$ for all $p$ in this neighborhood,
then $a+b$ is constant (because all $D\cap E_p$ are linearly
equivalent to $H$). On the other hand $a$ and $b$ can not be
constant. Hence $a-b$ is non-constant.
If we now define an equivalence relation $\sim$ on $E$
via $z\sim -z$, we obtain a globally well-defined holomorphic
non-constant map $E\to E/\!\sim$ locally given by $p\mapsto [a(p)-b(p)]$.
As a non-constant holomorphic map between compact Riemann surfaces
($E/\!\sim\,\simeq\P_1$),
this map must be surjective.
It follows that there exists a point $q\in E$ such that
$D\cap E_q=2\{r\}$ for some $r\in E_q$.
   Fix $F=E_q\subset A$ and define $F'=\pi^{-1}(F)\subset X$.
Then $F'\to F$ is ramified at exactly one point.
In local coordinates $x$ for the base and $t$ for the fiber of $L(D_0|_F)$
this implies $F'=\{(x,t):t^2=g(x)\}$ for some holomorphic
   function $g$ vanishing at the chosen base point $r$ of order $2$.
Then $g(x)=x^2e^{h(x)}$ for some holomorphic function $h$
and consequently $F'$ locally decomposes into two irreducible
components given by $t= \pm \, xe^{h(x)/2}$. For each of the two
components the projection onto $F$ is unramified.

Thus, if $\hat F$ denotes the normalization of $F'$ then
the naturally induced map $\hat F\to F$ is an unramified covering.
It follows that $\hat F$
is an elliptic curve, and therefore it follows that
$X$ contains an image of an elliptic curve under a non-constant
holomorphic map.
In particular, $X$ is not hyperbolic.
\end{proof}

\begin{thm}
\label{proj}
Let $E_i, 1 \leq i \leq q$ be smooth hypersurfaces of
the complex projective space $\pnc$ of dimension $n$
such that $E=\sum E_i$ is an s.n.c. divisor.
Assume that
\begin{enumerate}
\item
$q \geq n+1$.
\item
$\deg E \geq n+2$.
\end{enumerate}
Then every holomorphic curve $f: \C \to \pnc\setminus E$ is
algebraically degenerate.
\end{thm}

\begin{proof}
If $q \geq n+2$, then this is immediate from
Corollary~1.4. (iii) in \cite{nw02}.

Assume that $q=n+1$.
We observe that $K_\pnc+E$ is ample, because
$\deg(K_\pnc+E)\ge (-n-1)+n+2=1$. Thus $E$ is an s.n.c. divisor
for which $K_\pnc+E$ is ample, and therefore~$\pnc\setminus E$ is of log
general type.

Let $d$ be the l.c.m. of $d_i=\deg E_i, 1 \leq i \leq n+1$.
Define $k_i=d/d_i$.
Let $(x_0; \ldots ; x_n)$ be a homogeneous coordinate system
of $\pnc$.
Let $P_i(x_0, \ldots, x_n)$ be a homogeneous polynomial
defining $E_i$.

Now we define a morphism $\bar\pi$ from $\pnc$ to $\pnc$
as follows:
\begin{align*}
\bar \pi &: x \to (P_1^{k_1}(x); \ldots ;
  P_{n+1}^{k_{n+1}}(x)) 
\end{align*}

Then $\bar\pi$ is a finite morphism which restricts to a
finite morphism from $\pnc\setminus E$ to
$A=\{(y_0;\ldots ;y_n)\in\pnc: y_0\cdot y_1\cdot\ldots\cdot y_n \not= 0\}$.

Since $A$ is biholomorphic to $(\C^*)^{n}$, we may now use our
Main theorem and deduce that every holomorphic map 
$f: \C \to \pnc\setminus E$ is
algebraically degenerate.
\end{proof}

\begin{rmk}
\label{4.2}
{\rm
\begin{enumerate}
\item
Cf.\ \cite{nwy04} Theorem 7.1, in which $E_i, 1 \leq i \leq n$ are
hyperplanes and $\deg E_{n+1}\geq 2$.
\item
Grauert \cite{g89} dealt with the case where $n=2$ and
$E_i, 1 \leq i \leq 3$ are three smooth quadrics,
and proved that $\ptwoc \setminus E$ is Kobayashi hyperbolic.
The papers of Dethloff-Schumacher-P.M.\ Wong \cite{dswdj} and \cite{dswamj}
   followed it to make clear some part of the arguments in
Grauert \cite{g89}.
It was an essential step to prove the algebraic degeneracy of
$f: \C \to \ptwoc \setminus \sum_{i=1}^3 E_i$.
\item
Dethloff-Lu \cite{dl04} deals with the special case of $\dim X=2$ and
Brody curves $f: \C \to X$.
\end{enumerate}
}
\end{rmk}
{\it Examples. }  Let $D=\sum_{i=1}^3 D_i$
be the union of three quadrics of
$\ptwoc$ in sufficiently generic configuration
in the sense of \cite{dswamj} Proposition 4.1.
Then $X=\ptwoc\setminus D$ serves for an example of
Theorem \ref{hyp2} (cf.\ Remark \ref{4.2}).

\section{Strong Green-Griffiths conjecture}
In \cite{gg79} Green and Griffiths conjectured the following:

\begin{conjecture}
Let $X$ be a projective complex variety of general type. Then every
holomorphic map from $\C$ to $X$ is algebraically degenerate.
\end{conjecture}

This can be strengthened as follows:

\begin{conjecture}
Let $X$ be a projective complex variety of general type. Then there
exists a closed subvariety $E\subsetneq X$ such that $E$ contains
the image $f(\C)$ for every non-constant 
holomorphic map $f:\C\to X$.
\end{conjecture}

Here  we deal with this ``strong Green-Griffiths' conjecture'' for surfaces.

\begin{thm}
\label{stGG}
Let $A$ be a semi-abelian surface, and
let $X$ be a smooth surface of log general type.
Assume that there exists a proper finite
morphism $\pi:X \to A$.
Then there are only finitely many non-hyperbolic curves $C$ on $X$;
moreover an arbitrary nonconstant holomorphic curve
$f: \C \to X$ has the image contained in one of such $C$'s.
\end{thm}

\begin{proof}
The ramification locus $R$ of $\pi$ is the set of all points
$x\in X$ for which the differential $d\pi_x:T_x X \to T_{\pi(x)} A$
is not surjective and we denote by $R^*$
the set of all point $x\in X$ for which $d \pi_x$ has exactly
rank one.
Set $D=\pi(R)$.
In the same way as in \S3 (2) we see that $\std$ is trivial,
since $X$ is of log general type.
We choose a smooth equivariant compactification $A\hookrightarrow\bar A$
in which $D$ is in general position (\cite{nwy04} \S3).
Then we choose a compactification $\bar X$ of $X$
such that $\pi:X\to A$ extends to a finite morphism from $\bar X$
to $\bar A$ (see Proposition \ref{simple-cpt}).

Notice that there are only finitely many $A$-orbits in $\bar A$
(\cite{nwy04} Lemma 3.12).

Now we observe the following: If $E$ is a semi-abelian
subvariety of $A$ for which the fixed point set $\bar A^E$
in $\bar A$ is larger than $\bar A^A$,
then $E$ is the connected component
of the isotropy group of $A$ at a point in a one-dimensional
$A$-orbit in $\bar A$. Since $A$ is commutative, all points
in the same orbit have the same isotropy group.
Therefore the finiteness of the number of $A$-orbits implies
that there are only finitely many semi-abelian
subvarieties $E\subset A$ with $\bar A^E\ne\bar A^A$.

If $C$ is a non-hyperbolic curve on $X$,
then $C$ is either an elliptic curve or $\C^*$,
because a morphism from $\C$ or $\P_1$ to a semi-abelian variety
must be constant.
Hence the image of $C$ by $\pi$ is necessarily
an orbit in $A$ of a semi-abelian subvariety of $A$.
Thus it suffices to show that there exists only finitely many
such semi-abelian subvarieties of $A$ over which we can find
a non-hyperbolic curve on $X$.

Thus we have to investigate semi-abelian subvarieties $E\subset A$
of dimension one with an orbit $E(q) \subset A$ such that $\pi^{-1}(E(q))$
contains a non-hyperbolic curve.

We may therefore assume that $\bar A^E=\bar A^A$.

Now $E$ is either an elliptic curve or $\C^*$ and in both cases
$C\to E$ is an unramified covering. In particular $d \pi_x$
maps $T_xC$ surjectively on $T_{\pi(x)}E$ for every $x\in C$.
Therefore $C \cap R\subset R^*$.

We claim that furthermore $\bar C \cap \bar R \subset C\cap R$.
Indeed, if $C \ne\bar C$, then $\bar C=\P_1$ and $\pi(\bar C\setminus C)
\subset \partial A$. Therefore $\bar C\cap R=C\cap R$.
If $C=\bar C$, i.e., $C$ is an elliptic curve,
$\pi(\bar C\cap\partial X)\subset\bar A^E=\bar A^A$,
because $\pi(\bar C)$ is an $E$-orbit.
Since $D$ is generally positioned,
we have $\pi(\bar R)\cap \bar A^A=\emptyset$. Hence
$\bar C\cap\bar R\cap\partial X=\emptyset$.
Together, these arguments yield
\[
\bar C\cap\bar R\subset C\cap R
\subset C\cap R^*
\]

Next we define a ``Gauss map'' on $\overline{R^*}$:
We set
\[
\gamma: x \mapsto (\text{ Image of $d \pi_x$})\subset\P(\Lie A)
\]
It is readily verified that $\gamma$ is a rational map, thus it extends
to a morphism from the closure $R'=\overline{R^*}$ to $\P_1\simeq\P(\Lie A)$.
Let $R_1$ denote the union of irreducible components of $R'$
on which $\gamma$ is locally constant.
Then each irreducible component $K$ of $R_1$ maps onto an orbit
of a one-dimensional algebraic subgroup $H$ of $A$ in $\bar A$.
The value of the Gauss map is evidently $\Lie H$.
But then $x\in K\cap E$ implies that $E=H$ and that there is an $E$-orbit
inside $D$ which is excluded. Therefore $C\cap R_1=\emptyset$.
Let $R_2$ denote the union of all irreducible components of $R'$
along which $\gamma$ is nowhere locally constant.
Let $x\in C\cap R_2$. A calculation in local coordinates shows
that $\mult_x(C,R)=\mult_x\gamma+1$.
Therefore $\deg(C\cap R)\le 2\deg\gamma$.

In this way we obtain a universal bound for the degree of $\bar C$
with respect to the big divisor $\bar R$.
Therefore all such curves $C$ are contained
in a finite number of families. Now fix a one-dimensional semi-abelian
subvariety $E$ and consider all such curves $\bar C$ for which $\pi(C)$
is an orbit of the fixed algebraic subgroup $E$ of $A$.
Then $\pi(C)$ is an $E$-orbit containing
$\pi(R_2\cap C)$. Therefore $\pi(C)$ is uniquely determined by $C\cap R_2$
provided the latter is not empty.
Now $\{x\in R_2:\gamma(x)=[\Lie E]\}$ is finite and the number of
curves with empty intersection with $\bar R$ is finite as well,
since $\bar R$ is big.
Thus we have finished the proof that there are only finitely many such
curves.

The last assertion for entire holomorphic curves follows from
the Main Theorem combined with the above obtained result.
\end{proof}

\bigskip
\baselineskip=12pt
\rightline{J. Noguchi}
\rightline{Graduate School of Mathematical Sciences}
\rightline{University of Tokyo}
\rightline{Komaba, Meguro,Tokyo 153-8914}
\rightline{Japan}
\rightline{e-mail: noguchi@ms.u-tokyo.ac.jp}
\bigskip
\rightline{J. Winkelmann}
\rightline{Institut \'Elie Cartan}
\rightline{Universit\'e Nancy I}
\rightline{B.P.239}
\rightline{54506 Vand\oe uvre-les-Nancy Cedex}
\rightline{France}
\rightline{e-mail:jwinkel@member.ams.org}
\bigskip
\rightline{K. Yamanoi}
\rightline{Research Institute for Mathematical Sciences}
\rightline{Kyoto University}
\rightline{Oiwake-cho, Sakyoku, Kyoto 606-8502}
\rightline{Japan}
\rightline{e-mail: ya@kurims.kyoto-u.ac.jp}
\end{document}